\documentclass[12pt]{article}
\usepackage{a4}
\usepackage{amssymb}

\begin{document}
\newtheorem{theorem}{Theorem}[section]
\newtheorem{lemma}[theorem]{Lemma}
\newtheorem{corollary}[theorem]{Corollary}
\newtheorem{definition}[theorem]{Definition}
\newtheorem{proposition}[theorem]{Proposition}
\newtheorem{defprop}[theorem]{Definition-Proposition}
\newtheorem{example}[theorem]{Example}
\newtheorem{remark}[theorem]{Remark}
\newcommand{\Proof}{\noindent{\bf Proof:} }
\catcode`\@=11
\@addtoreset{equation}{section}
\catcode`\@=12
\renewcommand{\theequation}{\arabic{section}.\arabic{equation}}

\def\sqr#1#2{{\vcenter{\vbox{\hrule height.#2pt\hbox{\vrule width.#2pt
  height#1pt \kern#1pt \vrule width.#2pt}\hrule height.#2pt}}}}
\def\square{\mathchoice\sqr64\sqr64\sqr{2.1}3\sqr{1.5}3}
\def\argh{\rightharpoonup}
\def\m#1{m_{(#1)}}
\def\h#1{h_{(#1)}}
\def\k#1{k_{(#1)}}
\def\M{{\cal M}}
\def\F{{\cal F}}
\def\mm#1{m_{(#1*)}}
\def\a#1{a_{(#1)}}
\def\aa#1{a_{(#1*)}}
\def\b#1{b_{(#1)}}
\def\bb#1{b_{(#1*)}}
\def\l#1{l_{(#1)}}
\def\n#1{n_{(#1)}}
\def\nn#1{n_{(#1*)}}
\def\totimes{\tilde{\otimes}}
\def\c#1{c_{(#1)}}
\def\cc#1{c_{(#1*)}}
\def\id{{\rm id}}
\def\ydh{{\cal YD}(H)}
\def\yddh{{\cal YD}(H^*)}
\def\lefth{_{H}{\cal M}}
\def\lefthd{_{H^*}{\cal M}}
\def\righth{{\cal M}^H}
\def\righthd{{\cal M}^{H^*}}
\def\C{{\cal C}}
\def\D{{\cal D}}
\def\R{{\cal R}}

\title{SOME ISOMORPHISMS FOR THE BRAUER GROUPS OF A HOPF ALGEBRA}
\author{Giovanna Carnovale}
\date{\baselineskip 14pt Department of Mathematics and Computer Science\\
Universitaire Instelling Antwerpen\\
Universiteitsplein 1 - 2610 Wilrijk - Belgium\\
email: gcarno@uia.ua.ac.be\\}

\maketitle

\begin{abstract}
Using  equivalences 
of categories 
we provide isomorphisms between the Brauer groups of different Hopf algebras. 
As an example, we show that when $k$ is a field of 
characteristic different from $2$
 the Brauer groups $BC(k,\,H_4,\,r_t)$ for every dual 
quasitriangular structure $r_t$ on Sweedler's Hopf algebra $H_4$ are all 
isomorphic to the direct sum of $(k,+)$ and the Brauer-Wall group of $k$.
We  provide an 
isomorphism between
 the Brauer group of a Hopf algebra $H$ and the Brauer group of the dual Hopf algebra $H^*$
 generalizing a result of Tilborghs. Finally we relate the Brauer groups of $H$ and of its 
 opposite and co-opposite Hopf algebras.\\
\end{abstract}

\section{Introduction}

There have been given several generalizations of the Brauer group of a field $k$, due among others to 
Wall, Long, Van Oystaeyen, Caenepeel and Zhang. In particular, Caenepeel, Van Oystaeyen and Zhang
defined in \cite{CVZ} the Brauer group $BQ(k,\,H)$ of
 a Hopf algebra $H$ with bijective 
antipode. This is a special case of 
 Brauer group of a  braided monoidal category
(see \cite{VZ}: the Brauer group of a symmetric monoidal category 
had been defined by B. Pareigis in \cite{Par}). Here the category is that
of  left modules of the  Drinfel'd quantum double 
(see \cite{drin} and  \cite{Maj}) of a finitely generated projective
 Hopf algebra $H$ over a commutative ring $k$.\\
 $BQ(k, H)$ generalizes the Brauer-Long group of a 
commutative and cocommutative Hopf algebra over a commutative ring defined by Long in \cite{Long}. 
In fact it is shown in \cite{CVZ} that
when $H$ is a commutative cocommutative Hopf algebra, the Brauer group of $H$ and the Brauer-Long group of $H$ are
(anti-)isomorphic. \\
Another example of Brauer group of a braided monoidal category is the Brauer group $BC(k,\,H,\,r)$ 
where $H$ is a dual 
quasitriangular Hopf algebra with universal $R$-form $r$. In this case the category is that of  right $H$-comodules.\\
By results in \cite{doi} and results in \cite{mabook} (the dual version can be found in \cite{drin2} and the
survey book \cite{Kas} in the context of quasi
Hopf algebras)
 twisting the algebra structure of $H$ by a $2$-cocycle $\sigma$ provides a new dual quasitriangular
Hopf algebra $_{\sigma}\!H_{\sigma^{-1}}$ whose
comodule category is equivalent to that of $H$. Hence the two Brauer groups $BC(k,\,H,\,r)$ and 
$BC(k,\,_{\sigma}H_{\sigma^{-1}},\,r_\sigma)$ will be isomorphic. 
In this paper we use this result in order to show that although the universal $R$-forms $r_t$ for $t\in k$
 of Sweedler's four dimensional Hopf algebra $H_4$
are not all isomorphic (this fact was proved by Radford in \cite{Rad}), the Brauer groups $BC(k, H, r_t)$ are all
 isomorphic
for every $t\in k$.
This result is achieved by finding a suitable cocycle which does not change the algebra structure of $H_4$ but changes the 
universal $R$-form $r_t$ into $r_0$. Since the Brauer group $BC(k,H_4, r_0)$ (or equivalently, $BM(k, H_4, R_0)$ 
the Brauer group of the braided category of $H_4$-modules with braiding given by the triangular $R$-matrix $R_0$)
has been computed in \cite{VOZ3}, the computation of all $BC$'s for  $H_4$ is accomplished. 
This group is the direct sum of 
the additive group $(k,+)$ and of another generalization of the Brauer group of $k$, the so-called
 Brauer-Wall group $BW$ of $k$.  $BW(k)$ is the Brauer group of the category of ${\bf Z}_2$-graded $k$-modules and 
 it has been explicitely described by Wall in \cite{wall} by means of an exact sequence having the usual Brauer group of $k$ as a kernel.\\\\
Cocycle twisting also provides an isomorphism between the groups $BC(k,\,H,\,r)$ and $BC(k,\,H^{op},\,r\tau)$
 where 
$\tau$ is the usual flip
and $H^{op}$ denotes the Hopf algebra with opposite product.\\
In the case of a commutative and  cocommutative finitely generated and projective Hopf algebra $H$, 
it was shown by F. Tilborghs in \cite{Til} that the Brauer-Long group
of $H$ is (anti-)isomorphic to the Brauer-Long group of $H^*$.\\
Hence for $H$ commutative and cocommutative, finitely generated and projective one has an isomorphism between
the Brauer group $BQ$ of $H$ and the Brauer group $BQ$ of $H^*$. We shall use an isomorphism of Radford involving 
$D(H)$ and $D(H^*)$ where $D(H)$ denotes  the Drinfel'd quantum double of $H$,
together with  the results about the Brauer group $BC$ of the opposite Hopf algebra in order to generalize 
Tilborghs' result to the case of finitely generated  projective Hopf algebras with a bijective antipode. 
Everything boils down to  the fact that there is 
an anti-equivalence of braided categories between the category of Yetter-Drinfel'd $H$-modules 
(or crossed bimodules, or Quantum Yang-Baxter modules) and that of Yetter-Drinfel'd $H^*$-modules given by the usual
duality functor mapping left $H$-modules to right $H^*$-comodules and right $H$-comodules to left $H^*$-modules. 
Equivalently, 
there is an anti-equivalence of braided categories between left $D(H)$-modules and left $D(H^*)$-modules.\\
At the end of the paper we relate the Brauer groups of $H$ and $H^{op}$. We show that there is always an isomorphism
from $BQ(k,H)$ to $BM(k, D(H^{op}), \tau \R_{D(H^{op})}^{-1})$
where $\R_{D(H^{op})}$ denotes the standard $R$-matrix of $D(H^{op})$.
 If $D(H)$ were triangular with respect to its standard
 $R$-matrix,
an isomorphism $BQ(k,H)\simeq BQ(k,H^{op})$ would follow because $BQ(k,H)\simeq BM(k,D(H),\R)$. However
for nontrivial $H$
the quantum double  $D(H)$ cannot be triangular with respect to its standard $R$-matrix.\\
The paper is organized as follows. In Section \ref{category}  the construction of
 the Brauer group of a braided monoidal category is recalled. The particular case of the Brauer group(s) of a 
 Hopf algebra is included in a subsection.
 In Section  \ref{cocycle}  the theory of cocycle twists of a Hopf algebra is used 
 in order to get isomorphisms between  the groups $BC$ of Hopf algebras related by a twist. 
 The main results of the paper are to be found in Sections
 \ref{h4} and \ref{hop}. In Section \ref{h4} I apply the isomorphism described above to the particular example of  Sweedler's 
 four dimensional Hopf algebra $H_4$ in order to obtain $BC(k,H_4,r_t)$ for every $t\in k$.
  In Section \ref{hop}, I provide an isomorphism between  the Brauer group $BQ$ of a Hopf algebra $H$ with
 the Brauer group $BQ$ of its dual and I relate them to the Brauer groups of its opposite Hopf algebra.

\section{The Brauer group of a braided category}\label{category}

The Brauer group of a braided monoidal category was defined in \cite{VZ}
and this definition contains all known Brauer groups. The case of a symmetric category  had been already treated 
in \cite{Par} and in the symmetric case the Brauer group is abelian. Here we give a short account of the general
construction.\\
Let $\cal C$ denote a braided monoidal category with $\otimes$, $\psi$ and $I$ respectively the tensor,
braiding and identity object. For objects $P$ and $Q$ in $\cal C$ if the functor ${\cal C}(-\otimes P,\, Q)$
is representable one denotes the representing object by $[P,\,Q]$. Similarly if the functor ${\cal C}(P\otimes -,\,Q)$ 
is representable, 
the representing object is denoted by $\{P,\,Q\}$. 
An object $A$ in $\C$ is an algebra if there are morphisms $m\colon A\otimes A\to A$ (product)
 and $\eta\colon I\to A$ (unit)
satisfying associativity an unitary conditions. The $\C$-opposite algebra $\bar A$ is then defined as $A$ as object but
with product $\bar m:=m\circ\psi$ and same unit. The tensor product of two algebras
$A$ and $B$ in $\cal C$ becomes an algebra in $\cal C$ denoted by $A\# B$
with product $(m_A\otimes m_B)\circ(\id\otimes\psi\otimes\id)$. For more details, see \cite{mabook}
 and references therein. In particular, $[P,P]$ and $\{P,\,P\}$ are algebras in $\C$.\\
 For an algebra $A$ in $\C$ one has the two maps:
 \begin{eqnarray}F\colon \C(X,\, A\#{\bar A})&&\to \C(X,\, [A,\,A])\\
 F(a\# {\bar b})<d>&&=a{\bar m}(b\otimes d)\in\C(X\otimes Y,\,A)
 \end{eqnarray} for every $d\in\C(Y,\,A)$
 where $X$ and $Y$ are objects in $\C$ and $<d>$ stands for ``evaluation'' at $d$
and  
\begin{eqnarray}G\colon \C(X,\, {\bar A}\#{ A})&&\to \C(X,\, \{A,\,A\})\\
 <d>G({\bar a}\# b)&&={\bar m}(d\otimes a)b\in\C(X\otimes Y,\,A)
 \end{eqnarray} for every $d\in\C(Y,\,A)$
 where $X$ and $Y$ are objects in $\C$ and $<d>$ stands for ``evaluation'' at $d$. 
An algebra in $\C$ is called $\C$-Azumaya if $F$ and $G$ are isomorphisms and $A$ is faithfully projective in $\C$
(see \cite{Par} or \cite{VZ} for the definition of faithfully projective).
It turns out that: the product of two $\C$-Azumaya algebras is $\C$-Azumaya; the opposite algebra of a $\C$-Azumaya algebra
is $\C$-Azumaya. If $P$ is 
faithfully projective, then $[P,\,P]$ is also $\C$-Azumaya.
One can define an equivalence relation on the set of $\C$-Azumaya algebras:
$A\sim B$ if there exist faithfully projective objects $M$ and $N$ such that
\begin{equation}A\#[M,\,M]\simeq B\#[N,\,N].\end{equation}
It is proved in \cite{VZ} that this is indeed an equivalence relation and that the set of equivalence classes
becomes a group $Br(\C)$ with product induced by $\#$. The inverse of a class represented by an algebra $A$ will
be the class represented by the algebra $\bar A$.  \\
In \cite{VZ} the second Brauer group $Br'(\C)$   was also defined if $\C$ satisfies some extra 
conditions. Under  general conditions there is a group homomorphism $Br'(\C)\to Br(\C)$ 
which is very often a monomorphism and it is the identity if the unit object $I$ is projective.
The elements of $Br'(\C)$ are classes of  separable $\C$-Azumaya algebras in the category (see \cite{Par} or 
\cite{VZ}
for the definition).
\begin{remark}{\rm  It is clear that if we replace a category $\C$ by its opposite category $\C^{op}$
  with $A\otimes^{op} B:=B\otimes A$ for two objects in $\C^{op}$, $f\otimes^{op} g:=g\otimes f$ for two morphisms
  $f$ and $g$, and $\psi^{op}_{AB}:=\psi_{BA}$, then the map 
  $$\alpha\colon Br(\C)\to Br(\C^{op})$$
   $$[A]\mapsto [\bar A]_{op}$$
 where $[A]$ denotes the class of $A$ and $[A]_{op}$ denotes the equivalence class of $A$ in $\C^{op}$
  is a well-defined isomorphisms of groups which induces  also an isomorphism  between 
  $Br'(\C)$ and $Br'(\C^{op})$. In fact it is clear that $A$ is $\C$-Azumaya iff $\bar A$ is $\C^{op}$-Azumaya
  and $A\sim_{\C} B$ iff $\bar A\sim_{\C^{op}}\bar B$ iff $A\sim_{\C^{op}} B$.

}\end{remark}

\subsection{The Brauer groups of a Hopf algebra}

From now on  $k$ will denote a commutative ring.
Every $k$-module $M$ will be  assumed to be finitely generated, projective and faithful.
This implies the existence of ``dual bases''
$\{m_i\}\subset M$ and $\{m^*_j\}\subset M^*$ for which $\sum m^*_j(m)m_j=m$ for every $m\in M$. \\
All tensor product will be intended to be over $k$. \\
All Hopf algebras 
will  be assumed to be  finitely generated, projective and with bijective antipode $S$.  
$<-,\,->$ shall always denote evaluation between $H^*$ and 
$H$. Since $H$ is finitely generated and projective, $H^*$ is also a finitely generated projective Hopf algebra.\\
A $k$-module $M$ is called faithfully projective (progenerator) if there exist elements 
$m^*_j\in M^*$ and $m_j\in M$ for which $\sum m_j^*(m_j)=1$.\\
 $H^{op}$ shall denote the Hopf algebra with opposite product and $H^{cop}$ shall denote the Hopf algebra with 
 opposite 
 coproduct.\\
The Brauer group $BQ(k,\,H)$ of a Hopf algebra $H$ over $k$ is a particular case of Brauer group
of a braided monoidal category where the category is that
of Yetter-Drinfel'd $H$-modules.
\begin{definition} A Yetter-Drinfel'd $H$-module $M$ is a $k$-module
which is a left $H$-module and a right $H$-comodule satisfying the compatibility condition
\begin{equation}\sum h_{(1)}\cdot m_{(0)}\otimes h_{(2)}m_{(1)}=
\sum(h_{(2)}\cdot m)_{(0)}\otimes (h_{(2)}\cdot m)_{(1)}h_{(1)}\end{equation}
for every $h\in H$ and $m\in M$. In the formula $\cdot$ denotes the $H$-action on $M$, 
$\chi(m)= \sum m_{(0)}\otimes m_{(1)}\in M\otimes H$
 is the right comodule structure map and $\Delta(h)=\sum h_{(1)}\otimes h_{(2)}$ denotes the coproduct in $H$.
\end{definition}

\begin{remark} {\rm A Yetter-Drinfel'd module is also sometimes called a crossed bimodule (see \cite{Maj})
or a Quantum Yang-Baxter $H$-module (see \cite{LR}).}
\end{remark}
It is a result of S. Majid in \cite{Maj} that the category of Yetter-Drinfel'd $H$-modules is equivalent to the category of 
left $D(H)$-modules where $D(H)$ is the Drinfel'd double of $H$ defined in \cite{drin}.\\
In order to fix notation we recall that $D(H)$ is a quasitriangular Hopf algebra whose underlying coalgebra is $H^{*,cop}\otimes H$,
with product $$(\xi\otimes a)(\eta\otimes b)=\sum \xi\eta_{(2)}\otimes \a2 b <\eta_{(1)},\,S^{-1}(\a3)><\eta_{(3)},\,\a3>$$
and with $R$-matrix $\sum (\varepsilon\otimes h_i)\otimes (h^*_i\otimes 1)$ where $\{h_i\}$ and $\{h^*_i\}$ 
are dual bases in $H$ and $H^*$.\\
In this notation the $H$-action $\cdot$ and the $H^*$-action
$h^*\argh m=(\id\otimes<h^*,\,->)\chi$ together with the compatibility condition define a $D(H)$-module structure on 
a Yetter-Drinfel'd module $M$ and viceversa because $H$ and $H^*$ are subalgebras of $D(H)$.
 We shall denote the $D(H)$-action by $\triangleright$.\\
It is well-known that the tensor product of two Yetter-Drinfel'd $H$-modules $M$ and $N$ can be naturally equipped of 
a Yetter-Drinfel'd $H$-module
 structure denoted by $M\totimes N$ as follows:
\begin{equation}\chi (a\otimes b)=\sum \a0\otimes\b0\otimes \b1\a1\end{equation}
 and 
 \begin{equation}h\cdot(a\otimes b)=\sum\h1\cdot a\otimes \h2\cdot b.\end{equation}
 The category $\ydh$ of Yetter-Drinfel'd modules
   (with module and comodule morphisms) together with $\totimes$ becomes a monoidal category.
Moreover, there is an isomorphism of Yetter-Drinfel'd $H$-modules between $M\tilde\otimes N$ and 
$N\tilde\otimes N$  given by 
\begin{equation}\phi_{MN}(m\totimes n)=\sum \n0\totimes \n1\cdot m\end{equation}
 which makes of
the category of Yetter-Drinfel'd $H$-modules a braided monoidal category. 
The mentioned equivalence of ${\cal YD}(H)$ and $_{D(H)}{\cal M}$ is an equivalence of braided monoidal categories.
In terms of $D(H)$-modules, if
 $\R=\sum \R_1\otimes \R_2$ is
the standard $R$-matrix for $D(H)$ the braiding is nothing but
\begin{equation}\phi_{MN}(m\totimes n)=\sum \R_2\triangleright n\totimes \R_1\triangleright m\end{equation}
\\
In this setting for a Yetter-Drinfel'd module $P$, $[P,\,P]$ and $\{P,\,P\}$ are the usual
 $End(P)$ and $End(P)^{op}$ respectively, equipped 
with the Yetter-Drinfel'd module structures:
 \begin{equation}(h\cdot f)(m)=\sum \h1\cdot f(S(\h2\cdot m));\label{dritto1}\end{equation}
  \begin{equation}\chi(f)(m)=\sum f(\m0)_{(0)}\otimes S^{-1}(m_{(1)})f(\m0)_{(1)}\label{dritto2}\end{equation}
  for every $m\in M$ and $f\in End(M)$ for $End(M)$
  and 
  \begin{equation}(h\cdot' f)(m)=\sum \h2\cdot f(S^{-1}(\h1)\cdot m);\label{rovescio1}\end{equation}
  \begin{equation}\chi'(f)(m)=\sum f(\m0)_{(0)}\otimes f(\m0)_{(1)}S(\m1)\label{rovescio2}\end{equation}
  for every $m\in M$ and $f\in End(M)$ for $End(M)^{op}$. Those constructions come from two possible 
natural Yetter-Drinfel'd
   module structures on
   $[P,\,I]=P^*$.\\
An algebra in $\ydh$ is called a Yetter-Drinfel'd
 $H$-module
 algebra and it corresponds to a $D(H)$-module algebra:
\begin{definition} A Yetter-Drinfel'd $H$-module algebra $M$ is an algebra having the structure of a  Yetter-Drinfel'd
 $H$-module 
and such that the module and comodule structure make of $M$ a left $H$-module algebra and a right 
$H^{op}$-comodule algebra. 
\end{definition}
The tensor product $A\totimes B$ of two Yetter-Drinfel'd $H$-module algebras $A$ and $B$ becomes a Yetter-Drinfel'd module algebra
denoted by $A\# B$ with product given by  
 $(a\# c)(b\# d):=a\phi(c\otimes b)d=\sum a\b0\#(\b1\cdot c)d$.\\
The maps $F$ and $G$ become then:
 \begin{eqnarray}F\colon A\# {\bar A}&&\to End(A)\\
 F(a\# \bar b)(c)&&=\sum a\c0(\c1\cdot b)=m(a\otimes\phi_{AA}(b\otimes c))
 \end{eqnarray}
 and 
 \begin{eqnarray}G\colon \bar A\# { A}&&\to End(A)^{op}\\
 G(\bar a\# b)(c)&&=\sum \a0(\a1\cdot c)b=m(\phi_{AA}\tau\otimes\id)
 (\id\otimes\tau)(a\otimes b\otimes c)
 \end{eqnarray}
 for every $a,\,b$ and $c$ in $A$. It had already been proved in \cite{CVZ} that $F$ and $G$ are Yetter-Drinfel'd $H$-module
  algebra maps.
  A Yetter-Drinfel'd $H$-module algebra $A$ is  $H$-Azumaya if $F$ and $G$ are isomorphisms. 
 In this case $A\sim B$ if there exist faithfully projective
 Yetter-Drinfel'd $H$-modules $M$ and $N$ such that $A\# End(M)\simeq B\# End(N)$ as Yetter-Drinfel'd
  $H$-module algebras. $Br(\ydh)$ in usually denoted by $BQ(k,\,H)$. \\
  An $H$-Azumaya algebra is said to be strongly $H$-Azumaya if $k$ is a direct summand of $A$ as Yetter-Drinfel'd $H$-module. 
  The second Brauer group $Br'(\ydh)$ is usually denoted $BQS(k,\,H)$ and it coincides with $BQ(k,\,H)$ 
  if $H$ is semisimple-like and cosemisimple-like (see Prop. 2.27 in \cite{CVZ2}). In general $BQS(k,\,H)$
   is the subgroup whose elements are classes parametrized by strongly $H$-Azumaya algebras.\\ 
  
  If $H$ is  a quasitriangular Hopf algebra with universal $R$ matrix
$R=\sum R_1\otimes R_2$
 it is well known that to every module (algebra) $A$ one can associate a right $H^{op}$-comodule 
(algebra) structure on $A$ given by
\begin{equation}\label{comodule}\chi_R(a)=\sum R_2\cdot a\otimes R_1\end{equation}  
obtaining a Yetter-Drinfel'd $H$-module algebra. 
In this case the category $_H\!{\cal M}$ of left $H$-modules with $H$-module maps
 is a full subcategory of
 $\ydh$. 
$Br(_H{\cal M})$ is then a subgroup of $BQ(k,\,H)$ and it is usually denoted by $BM(k,\,H,\,R)$.
It is  the subgroup of $BQ(k,\,H)$ whose elements are represented by a Yetter-Drinfel'd module
 algebras whose comodule structure is defined by (\ref{comodule}).\\ Dually,
 if $H$ is a dual quasitriangular Hopf algebra with universal $R$-form $r$,
 to a right $H^{op}$-module (algebra) $A$ one can associate a
left $H$-module (algebra) structure on $A$ given by 
\begin{equation}\label{module}h\cdot a=\sum\a0 r(h\otimes \a1)\end{equation} obtaining a 
Yetter-Drinfel'd module algebra. The category 
${\cal M}^{H}$ of right $H$-comodules  with $H$-comodule maps is a full subcategory of $\ydh$. $Br({\cal M}^H)$ 
is then a subgroup of 
$BQ(k,\,H)$ denoted by $BC(k,\,H,\,r)$. It is the subgroup of $BQ(k,\,H)$ whose elements are represented by a 
Yetter-Drinfel'd module
 algebras whose module structure is defined by (\ref{module}).\\
It is well-known that $BQ(k,\,H)\simeq BC(k,\,D(H)^*,\, r)\simeq BM(k,\,D(H), \R)$ where $D(H)$ is the
Drinfel'd double of $H$ and $\R$ is its standard $R$-matrix. Hence it is enough to study $BC(k,\,H,\,r)$
 for a dual quasitriangular Hopf algebra.

\section{An equivalence of categories}\label{cocycle}

 In this section we show a few isomorphisms for the group $BC$ of a dual quasitriangular Hopf algebra.
  Everything can be dualized,
considering quasitriangular Hopf algebras and $BM$. We leave this task to the reader.\\

Let $H$ be a bialgebra over $k$ and let $B$ be a left (resp. right) $H$-comodule algebra with comodule map $\chi$. 
A left (resp. right) $2$-cocycle $\sigma$ is a linear map $\sigma\colon H\otimes H\to k$ satisfying 
\begin{itemize}
\item $\sum\sigma(k_{(1)}\otimes \m1)\sigma(h\otimes \k2\m2)=\sum\sigma(\h1\otimes \k1)
\sigma(\h2\k2\otimes m)$\\
$${\hbox{(resp. }}\quad
\sum\sigma(\k2\otimes \m2)\sigma(h\otimes \k1\m1)=\sum\sigma(\h2\otimes
 \k2)\sigma(\h1\k1\otimes m) )$$ 
\item $\sigma(h\otimes 1)=\sigma(1\otimes h)=\varepsilon(h)$  $\forall$ $h,\,k,\,m\,\in H$ 
\end{itemize}

Then, the $\sigma$-left (resp. $\sigma$-right) twisted comodule $_{\sigma}B$ (resp. $B_{\sigma}$) is an algebra with  
the same underlying vector space as $B$, and product given by:

 \begin{equation}\bar a\cdot\bar b=\sum\sigma(\a1\otimes \b1)\overline{\a0\b0}\end{equation}  
if $\chi(a)=\sum \a1\otimes \a0$, and $\chi(b)=\sum \b1\otimes \b0\in H\otimes B$ for $a,\,b\,\in\,B$;\\ 
  \begin{equation}{\rm (resp. }\quad\bar a\cdot\bar b=\sum\overline{\a0\b0}\sigma(\a1\otimes \b1)\end{equation} 
 if $\chi(a)=\sum \a0\otimes \a1$, and $\chi(b)=\sum \b0\otimes \b1\in B\otimes H$),  
 where  $a\mapsto\bar a$ denotes the identification of vector spaces (see for instance \cite{Mon}, Paragraph 7.5).\\
If $B$ is a bialgebra and $\sigma$ is a left 2-cocycle one can perform such a twist to $B$,
 viewed as a left (resp. right) $B$-comodule algebra. If $\sigma$ is a convolution invertible 
  left 2-cocycle, 
 then  $\sigma^{-1}$ is a right 2-cocycle. It is well-known (see for instance \cite{mabook},
\cite{doi} or \cite{doitake}) that the
  double twist $_{\sigma}B_{\sigma^{-1}}$ with the same coproduct of $B$ is again a bialgebra
   and if $B$ is a Hopf algebra, then $_{\sigma}B_{\sigma^{-1}}$
  is also a Hopf algebra with antipode $S_{\sigma}$ given by
 $(u\otimes S\otimes u^{-1})(\Delta\otimes\id)\Delta$. Here $u\in \,_{\sigma}\!B^*_{\sigma^{-1}}\simeq B^*$ 
 is the linear functional given by $u=\sigma(\id\otimes S)\Delta$.\\
 The following facts are also well-known
  (see for instance \cite{mabook} and references therein or, for a dual version involving the case of
   Drinfel'd
  quasi Hopf algebras, \cite{Kas} and references therein):
  \begin{itemize}
  \item The category of 
  right $H$-comodules is equivalent to the category of right $_{\sigma}H_{\sigma^{-1}}$-comodules as 
  monoidal category. In this case though we need the comodule structure of a tensor product to be
  \begin{equation}m\otimes n\to\sum \m0\otimes\n0\otimes\m1\n1.\end{equation}
   The monoidal functor $\cal F$ is the identity on objects
   (the coproduct is unchanged) and the natural transformation  
  $\eta\colon\F(M)\otimes\F(N)\to\F(M\otimes N)$ 
  \begin{equation}\eta(m\otimes n)=\sum \m0\otimes \n0 \sigma^{-1}(\m1\otimes \n1).\end{equation}
  The compatibility conditions for $\eta$ follow by the cocycle condition on $\sigma$.
  \item If $H$ is dual quasitriangular with universal $R$-form 
   $r$,
   then  $_{\sigma}H_{\sigma^{-1}}$ is also dual quasitriangular with 
   universal $R$-form given by $r_{\sigma}=(\sigma\tau)*r*\sigma^{-1}$.
  \item The categories $\M^{H}$ and 
   $\M^{_{\sigma}H_{\sigma^{-1}}}$
   are equivalent braided monoidal categories. In this case one uses the braiding
   \begin{equation}\phi_{MN}(m\otimes n)\mapsto\sum \n0\otimes\m0 r(\m1\otimes\n1)\end{equation}
    compatible with the different definition of 
   comodule structure on $M\otimes N$.
    One can check that the braiding
    $\phi_{MN}$ is respected by the functor $\cal F$ 
   and the natural transformation $\eta$ mentioned above.
   \item The universal $R$-forms $r$ and $r^{-1}\tau$ of a dual quasitriangular Hopf algebra
   $H$ are 2-cocycles and $_rH_{r^{-1}}=H^{op}=\, _{r^{-1}\tau}H_{r\tau}$. The new universal $R$-form
    is in both
   cases $r\tau$.
  \end{itemize}
  
  It is a straightforward check that if $H$ is a bialgebra with coproduct $\Delta$, and if $\sigma$ is 
  a 2-cocycle for 
  $H$, then $\sigma\tau$ is a 2-cocycle for $H^{op}$ and \\
   $(_{\sigma}H_{\sigma^{-1}})^{op}\simeq _{\sigma\tau}\!
  (H^{op})_{\sigma^{-1}\tau}$. Hence, repeating the discussion above for
   $H^{op}$ the category of right $H$-comodules with tensor product structure $$m\otimes n\mapsto 
  \sum \m0\otimes\n0\otimes\n1\m1$$ and braiding $$m\otimes n\mapsto\sum\n0\otimes\m0\,r(\n1\otimes\m1)$$
   is equivalent, 
  as braided monoidal category, to that of $_{\sigma}H_{\sigma^{-1}}$-comodules with tensor product 
  structure
  $$m\otimes n\mapsto\sum\m0\otimes\n0\otimes \sigma(\n1\otimes\m1)\n2\m2\sigma^{-1}(\n3\otimes\m3)$$ 
  and braiding
  $$m\otimes n\mapsto\sum\n0\otimes\m0\,\sigma(\m1\otimes\n1)r(\n2\otimes\m2)\sigma^{-1}(\n3\otimes\m3).$$ 
  The natural transformation $\eta$ is in this case: 
  \begin{equation}m\otimes n\mapsto \sum \m0\otimes\n0 \sigma^{-1}(\n1\otimes\m1).\end{equation}
  An algebra $A$ in $\M^{H}$ is mapped by the functor to the algebra $A_{\sigma^{-1}\tau}$ with product 
  $m_A\circ\eta_{A,A}$
  i.e. $a\cdot b:=\sum \a0\b0\sigma^{-1}(\b1\otimes\a1)$. Hence we have:
   \begin{proposition}\label{iso} Let $H$ be a (faithfully projective)
    dual quasitriangular Hopf algebra with universal $R$-form $r$. Let $\sigma$ be an invertible $2$-cocycle. \\
    Then
    $BC(k,\,H,\,r)\simeq BC(k,\, _{\sigma}H_{\sigma^{-1}},\,r_\sigma)$. The class of an $H$-Azumaya algebra 
    $A$ is mapped to the class of 
    $A_{\sigma^{-1}\tau}$. In particular, 
    $BC(k,\,H,\,r)\simeq BC(k,\, H^{op},\,r\tau)$.
    \end{proposition}
    \Proof It follows by the above observations. Observe that twisting $H$ by $r^{-1}\tau$ implies
    that $[A]$ is mapped to $[\overline{A^{op}}]$.\hfill$\square$
    \begin{remark}\label{gio}{\rm The result about $H^{op}$ could be obtained also by checking that
     the categories of right $H^{op}$
    comodules and of right $H$ comodules are anti-equivalent braided categories. Then the
     anti-isomorphism between the 
    two Brauer groups is given
    on representatives by $A\mapsto A^{op}$.}\end{remark}

    \section{An example: $BC(k,\,H_4,\,r_t)$}\label{h4}
     Let $k$ be a field of characteristic different from $2$ and let 
     $H_4$ denote Sweedler's four dimensional Hopf algebra over $k$ generated by $g$ and $h$ such that
     $g^2=1$, $h^2=0$ and $gh+hg=0$. As far as the coproduct $\Delta$ is concerned, $g$ is grouplike and $h$ is 
     twisted-primitive with
     $\Delta(h)=h\otimes g+1\otimes h$. The antipode $S$ is such that $S(g)=g$ and $S(h)=gh$.
     It is well-known that $H^*_4$ is isomorphic to $H_4$. An isomorphism is obtained sending $g$ to $f_1-f_g$
      and $h$ 
     to 
     $f_h+f_{gh}$ where $f_x$ denotes the dual element of $x\in H_4$. It is also well-known that $H_4$
      has a
      family of universal $R$-forms (and, dually of universal $R$-matrices) parametrized by the elements in $k$.
      The universal $R$-forms $r_t$ ($t\in k$) were firstly found by Radford in \cite{Rad} and
       are determined by the axioms of an $R$-form together with
      $$r_t(1\otimes x)=r_t(x\otimes 1)=\varepsilon(x)$$ 
      $$r_t(g\otimes g)=-1\quad r_t(g\otimes h)=r_t(h\otimes g)=r_t(g\otimes gh)=r_t(gh\otimes g)=0$$
      $$r_t(gh\otimes h)=r_t(h\otimes h)=r_t(gh\otimes gh)=-r_t(h\otimes gh)=t.$$
      A first observation is that all those structures are cotriangular, i.e. 
      $r_t*(r_t\tau)=\varepsilon\otimes \varepsilon=(r_t\tau)* r_t$.
   This means that for every $a$ and $b$ in $H_4$ one has:
   \begin{equation}\sum r_t(\a1\otimes\b1)r_t(\b2\otimes\a1)=\varepsilon(a)\varepsilon(b)=\sum r_t(\b1\otimes\a1)
   r_t(\a2\otimes\b2)\label{cotri}\end{equation}
   By the symmetry in the formula if we interchange $a$ and $b$ it is enough to check the left hand side
    condition on pairs of basis elements. We have various cases depending on the coproduct of $a$ and $b$:
   \begin{itemize}
   \item $a$ and $b$ are both grouplike elements of the basis (i.e. $1$ or $g$).  The left hand side of 
   $(\ref{cotri})$
    reads
     $r_t(a \otimes b)r_t(b\otimes a)$. This is $1^2$ for the pairs $(1,g)$, $(g,1)$ and $(1,1)$
     and it is $(-1)^2$ when $a=b=g$. In all cases this is equal to $\varepsilon(a)\varepsilon(b)$.
   \item One element in $\{a,\,b\}$ is grouplike and the other is twisted primitive (i.e. $h$ or $gh$).
   Then each summand of the left hand side of $(\ref{cotri})$
   will contain an expression of type $r_t(x\otimes y)$ with $x$ grouplike and $y$ 
   twisted primitive.  Therefore each summand is zero.
    For instance  $$r_t*r_t\tau(h\otimes g)=r_t(h\otimes g)r_t(g\otimes g)+r_t(1\otimes g)r_t(g\otimes h)=
   0=\varepsilon(g)\varepsilon(h).$$ Hence 
   $\sum r_t(\a1\otimes\b1)r_t(\b2\otimes\a1)=0=\varepsilon(a)\varepsilon(b)$.
   \item $a$ and $b$ are both twisted primitives.
   Then the only nonzero terms in the sum 
   on the left hand side of $(\ref{cotri})$ appear when  $\a1$ and $\b1$ are both twisted primitives, or when
    $\a1$ and $\b1$ are both grouplikes. If $a=b=h$ the expression becomes
   $$r_t*r_t\tau(h\otimes h)=r_t(h\otimes h)r_t(g\otimes g)+r_t(1\otimes 1)r_t(h\otimes h)=
   -t+t=\varepsilon(h)^2.$$
   If $a=b=gh$ the expression becomes
   $$r_t(gh\otimes gh)r_t(1\otimes 1)+r_t(g\otimes g)r_t(gh\otimes gh)=
   -t+t=\varepsilon(gh)\varepsilon(gh).$$
   If $a=h$ and $b=gh$ the expression becomes
   $$r_t(h\otimes gh)r_t(1\otimes g)+r_t(1\otimes g)r_t(gh\otimes h)=
   t-t=\varepsilon(gh)\varepsilon(h).$$
   Finally if $a=gh$ and $b=h$ the expression is
   $$r_t*r_t\tau(gh\otimes h)=r_t(gh\otimes h)r_t(g\otimes 1)+r_t(g\otimes 1)r_t(h\otimes gh)=
   t-t=\varepsilon(h)\varepsilon(gh).$$
  \end{itemize}   Hence $H_4$ is cotriangular for every universal $R$-form $r_t$. Therefore
  $BC(k,H_4,r_t)\simeq BC(k,H^{op},r_t^{-1})\simeq BC(k,H,r_t^{-1}\tau)$ and it is an abelian 
  group.\\
      Dually, one can check that $(H_4,\,R_t)$ is triangular for every 
     universal $R$-matrix $R_t$. Although triangularity of $R_t$
      follows by the previous result we sketch the proof for sake of completeness  because we have not
      met this result in the literature before. The family of $R$-matrices (see \cite{Rad} or \cite{LR}) 
      is given by
     $$R_t={1\over2}(1\otimes1+1\otimes g+g\otimes 1-g\otimes g)+{t\over2}
     (h\otimes h+h\otimes gh+gh\otimes gh-gh\otimes h)$$ for $t\in k$.
    $R_0=\tau R_0$ and $R_0^2=1\otimes 1$ because it corresponds to a Hopf involution 
    $f_R\colon H^*_4\to H_4^{cop}$ where $f_R(\xi)=(<\xi,->\otimes \id)(R_0)$ (see \cite{LR}). 
     Hence $(H_4,\,R_0)$ is triangular. Put $R_t=R_0+R'_t$. $(H_4,R_t)$ is triangular if
     $(\tau R_t)R_t=1\otimes 1=R_t(\tau R_t)$.
     Since $\tau$ is an algebra isomorphism it is enough to check the relation for $(\tau R_t)R_t$. 
     This expression is
     equal to
     $$(\tau R_0)R_0+(\tau R_0)R'_t+ (\tau R'_t)R'_0+(\tau R'_t)R'_t=1\otimes 1+(\tau R_0)R'_t+ (\tau R'_t)R'_0.$$
$(\tau R'_t)R'_t =0$  because $h^2=0$ appears in every component.
  Checking that  $(\tau R_0)R'_t+ (\tau R'_t)R'_0=0$ is a striaghtforward computation that we leave to the reader.\\ 
  The group $BM(k,\,H_4,\,R_0)$ has been
   computed by
  Van Oystaeyen and Zhang in \cite{VOZ3}. This group is isomorphic to  $BC(k,\,H_4,\,r_0)$ because the 
  universal $R$-matrix $R_0$ goes over to the universal $R$-form $r_0$ under the isomorphism 
  $H_4\to H_4^*$ previously given.\\
  We want to show here that $BC(k,\,H_4,\,r_t)\simeq BC(k,\,H,\,r_0)$ for every $t\in k$, hence that 
  $BM(k,\,H_4,\,R_0)\simeq
  BM(k,\,H_4,\,R_s)$ for every $s\in k$. We shall do this by providing for every $t\in k$ there exists
  a suitable element $\sigma_t\in (H_4\otimes H_4)^*$ such that 
  \begin{itemize}
  \item $\sigma_t$ is a left 2-cocycle for $H_4$;
  \item $\sigma_t$ is invertible;
  \item the twisted product in $_{\sigma_t}(H_4)_{\sigma_t^{-1}}$ coincides with the product in $H_4$;
  \item $\sigma_t\tau* r_t*\sigma_t^{-1}=r_0$.
  \end{itemize}
  The functional $\sigma_t$ is defined on the basis elements of $H_4\otimes H_4$ as follows:
  $$\sigma_t(x\otimes 1)=\sigma_t(1\otimes x)=\varepsilon(x)\quad{\hbox{for every $x\in H_4$.}}$$
  $$\sigma_t(g\otimes g)=1$$
  $$\sigma_t(h\otimes h)=\sigma_t(gh\otimes h)=-\sigma_t(h\otimes gh)=-\sigma_t(gh\otimes gh)={t\over2}$$
  and $\sigma_t(x\otimes y)=0$ whenever $x$ is grouplike and $y$ twisted primitive or the other way around.\\
 It is a $2$-cocycle if 
 for every triple of basis elements  $k,\,a,\,m$ there holds:
 $$\sum\sigma_t(k_{(1)}\otimes \m1)\sigma_t(a\otimes \k2\m2)=\sum\sigma_t(\a1\otimes \k1)
\sigma_t(\a2\k2\otimes m).$$
If one of the elements is $1$ then the condition is verified because $\sigma_t$ is unitary 
(i.e. it coincides with 
$\varepsilon$
on $1\otimes x$ and $x\otimes 1$). Hence we have to check the condition on all triples of elements in
 $\{g,\,h,\,gh\}$.\\
Then we have different cases depending on how often $g$ appears in the triple.
\begin{enumerate}
\item $g$ appears $3$ times.\\
Then we have $$\sigma_t(g\otimes g)\sigma_t(g\otimes 1)=1=\sigma_t(g\otimes g)\sigma_t(1\otimes g)$$
\item $g$ appears twice in the triple:\\
If $g=a=m$ and $k$ is twisted primitive the condition becomes
$$\sum\sigma_t(k_{(1)}\otimes g)\sigma_t(g\otimes \k2 g)=\sum\sigma_t(g\otimes \k1)
\sigma_t(g\k2\otimes g).$$
Since $k$ is twisted primitive, in every summand $\k1$ and $\k2$ can never be both grouplikes, 
hence both sums are $0$.
The cases $g=a=k$ and $g=k=m$ are checked similarly.
\item $g$ appears once in the triple:\\
If $g=m$ then the condition reads
$$\sum\sigma_t(k_{(1)}\otimes g)\sigma_t(a\otimes \k2 g)=\sum\sigma_t(\a1\otimes \k1)
\sigma_t(\a2\k2\otimes g).$$
The only nonzero component of the left hand side appears when $\k1$ is grouplike and $\k2=k$,
 hence the left hand side is equal to $\sigma_t(a\otimes kg)$. The only nonzero component of the right 
 hand side is 
 when $\a2\k2$ is grouplike, i.e. when both $\a2$ and $\k2$ are grouplikes, hence $\a1=a$, $\k1=k$
  and the right hand side
  becomes
 $\sigma_t(a\otimes k)$. It is straightforward to check that for the twisted primitives $a$ and $k$ in $H_4$
 there holds: $\sigma_t(a\otimes kg)=\sigma_t(a\otimes k)$.\\
 If $g=a$ by similar computations the left hand side becomes $\sigma_t(k\otimes m)$, with 
 $k$ and $m$ twisted primitives.
 The right hand side becomes $\sigma_t(gk\otimes m)=\sigma_t(k\otimes m)$ for $k$ and $m$ 
 twisted primitives.\\
 If $g=k$ the right hand side is $\sigma_t(a\otimes gm)$ and the left hand side is
  $\sigma_t(ag\otimes m)$. Again they 
 coincide
 for $a$ and $m$ twisted primitives.\\
 \item $g$ does not appear in the triple.\\
 Then $\sum\sigma_t(k_{(1)}\otimes \m1)\sigma_t(a\otimes \k2\m2)=0$ because if $\k2$ and $\m2$ are both 
 twisted primitives,
  their 
 product involves an $h^2$ hence it is zero, if one is twisted primitive and the other grouplike, then
  $\sigma_t(\k1\otimes\m1)=0$
 and if they are both grouplikes $\sigma_t(a\otimes\k2\m2)=0$. Similarly one shows that the 
 right hand side is also
  equal to
 zero.
 \end{enumerate}
Hence $\sigma_t$ is a left $2$-cocycle. Similarly one can prove that $\sigma_t$ is also a right $2$-cocycle.

We claim that $\nu_t\in(H_4\otimes H_4)^*$ is the convolution inverse of $\sigma_t$, where $\nu_t$
is defined on the basis elements as follows:
$$\nu_t(x\otimes 1)=\nu_t(1\otimes x)=\varepsilon(x)\quad{\hbox{for every $x\in H_4$.}}$$
  $$\nu_t(g\otimes g)=1$$
  $$\nu_t(h\otimes h)=\nu_t(gh\otimes h)=-\nu_t(h\otimes gh)=-\nu_t(gh\otimes gh)=-{t\over2}$$
  and $\nu(x\otimes y)=0$ whenever $x$ is grouplike and $y$ is twisted primitive or the other way around.\\
We show that $\sigma_t*\nu_t=\varepsilon\otimes\varepsilon=\nu_t*\sigma_t$. This means that we have to
 show that for every pair of basis elements  $a$ and $b$ in $H_4$ one has
 $$\sum \sigma_t(\a1\otimes\b1)\nu_t(\a2\otimes\b2)=\varepsilon(a)\varepsilon(b)
 =\sum\nu_t(\a1\otimes\b1)\sigma_t(\a2\otimes\b2).$$
Again we divide the different cases:
\begin{itemize}
\item $a$ and $b$ are both grouplikes i.e. $a,\,b\in\{1,\,g\}$. Then the expression becomes
$$\sigma_t(a\otimes b)\nu_t(a\otimes b)=1\cdot 1=\varepsilon(a)\varepsilon(b)=\nu_t(a\otimes b)\sigma_t(a\otimes b)$$
\item $a$ and $b$ are one grouplike and the other twisted primitive.
Then the expressions are always zero because in each summand there will be either a $\sigma_t(x\otimes y)$ or a
$\nu_t(x\otimes y)$ with one element grouplike and the other twisted primitive.
Hence $$\sum \sigma_t(\a1\otimes\b1)\nu_t(\a2\otimes\b2)=0=\varepsilon(a)\varepsilon(b)
 =\sum\nu_t(\a1\otimes\b1)\sigma_t(\a2\otimes\b2)$$
 \item $a$ and $b$ are both twisted primitive. The only nonzero terms in the sums will be
 those where both $\a1$ and $\b1$ are both grouplikes or both twisted primitives. We have
 for $a=h=b$:
  $$\sigma_t*\nu_t(h\otimes h)=\sigma_t(1\otimes 1)\nu_t(h\otimes h)+\sigma_t(h\otimes h)\nu_t(g\otimes g)+0=0=
  \varepsilon(h)^2$$
$$\nu_t*\sigma_t(h\otimes h)=\nu_t(1\otimes 1)\sigma_t(h\otimes h)+\nu_t(h\otimes h)
\sigma_t(g\otimes g)=0=\varepsilon(h)^2.$$
 For $a=h$ and $b=gh$:
 $$\sigma_t*\nu_t(h\otimes gh)=\sigma_t(1\otimes g)\nu_t(h\otimes gh)+\sigma_t(h\otimes gh)\nu_t(g\otimes 1)=0=
 \varepsilon(h)\varepsilon(gh)$$
 $$\nu_t*\sigma_t(h\otimes gh)=\nu_t(1\otimes g)\sigma_t(h\otimes gh)+\nu_t(h\otimes gh)\sigma_t(g\otimes 1)=0=
 \varepsilon(h)\varepsilon(gh).$$
 For $a=gh$ and $b=h$:
$$\sigma_t*\nu_t(gh\otimes h)=\sigma_t(gh\otimes h)\nu_t(1\otimes g)+\sigma_t(1\otimes g)\nu_t(gh\otimes h)
=0=\varepsilon(gh)\varepsilon(h)$$
$$\nu_t*\sigma_t(gh\otimes h)=\nu_t(gh\otimes h)\sigma_t(1\otimes g)+\nu_t(1\otimes g)\sigma_t(gh\otimes h)
=0=\varepsilon(gh)\varepsilon(h).$$
Finally for $a=b=gh$  
 $$\sigma_t*\nu_t(gh\otimes gh)=\sigma_t(g\otimes g)\nu_t(gh\otimes gh)+\sigma_t(gh\otimes gh)\nu_t(1\otimes 1)=0=\varepsilon(gh)^2$$
$$\nu_t*\sigma_t(gh\otimes gh)=\nu_t(g\otimes g)\sigma_t(gh\otimes gh)+\nu_t(gh\otimes gh)\sigma_t(1\otimes 1)=0=\varepsilon(gh)^2$$
\end{itemize}
 So $\nu_t=\sigma^{-1}_t$.
 Hence it makes sense to compute the product in $_{\sigma_t}H_{\sigma_t^{-1}}$. We shall see that 
 the product in  $_{\sigma_t}H_{\sigma_t^{-1}}$ coincides with the product in $H_4$.
  We check this fact
  on  products of the the generators $h$ and $g$: the other products follow by associativity.
 $$\bar g\cdot \bar g=\sigma_t(g\otimes g)\,\overline{g^2}\sigma_t^{-1}(g\otimes g)=1;$$
 $$\bar g\cdot\bar h=\sigma_t(g\otimes1)\,\overline{g}\,\sigma^{-1}_t(g\otimes h)+
 \sigma_t(g\otimes 1)\,\overline{gh}\,
 \sigma^{-1}_t(g\otimes g)+$$
 $$\sigma_t(g\otimes h)\,\overline{g^2}\,\sigma^{-1}_t(g\otimes g)=\overline{gh};$$
 $$\bar h\cdot\bar g=0+\sigma_t(1\otimes g)\,\overline{hg}\,\sigma^{-1}_t(g\otimes
  g)=\overline{hg};$$
 $$\bar h\cdot\bar h=
 \sigma_t(1\otimes 1)\,\overline{1}\,\sigma^{-1}_t(h\otimes h)+
 \sigma_t(1\otimes 1)\,\overline{h}\,\sigma^{-1}_t(h\otimes g)+$$
 $$ \sigma_t(1\otimes h)\,\overline{g}\,\sigma_t(h\otimes g)+
  \sigma_t(1\otimes 1)\,\overline{h}\,\sigma^{-1}_t(g\otimes h)+
   \sigma_t(1\otimes 1)\overline{h^2}\sigma^{-1}_t(g\otimes g)=$$   
$$\sigma_t(1\otimes h)\,\overline{hg}\,\sigma^{-1}_t(g\otimes g)+
\sigma_t(h\otimes 1)\,\overline{g}\,\sigma^{-1}_t(g\otimes h)+
\sigma_t(h\otimes 1)\,\overline{gh}\,\sigma^{-1}_t(g\otimes g)+$$
 $$ \sigma_t(h\otimes h)\,\overline{g^2}\,\sigma^{-1}_t(g\otimes g)=
 -{t\over2}+0+{t\over2}=0.$$
  Hence the product in  $_{\sigma_t}(H_4)_{\sigma_t^{-1}}$ coincides with the product in  $H_4$.\\
  Now we
   compute how $r_s$ changes under twisting, i.e. $\sigma_t\tau*r_s*\sigma^{-1}_t$ for every $s$ and $t$ in $k$.
   We shall prove that this
  is equal to $r_{t-s}$.\\
  It is known that $\sigma_t\tau*r_s*\sigma^{-1}_t$ is a universal $R$-form for the twist of $H_4$, 
  hence for $H_4$.
   By the properties of an $R$-form:
  $$r(a\otimes 1)=r(1\otimes a)=\varepsilon(a)\quad{\hbox{and}}$$
  $$r(ab\otimes c)=\sum r(a\otimes \c1)r(b\otimes \c2)\quad{\hbox{for every $a,\,b$ and $c\in H_4$}}$$ 
 it is enough to check  the equality when the first argument is $h$ or $g$. 
 Moreover, since $\Delta(h)$ and
  $\Delta(g)$
 can be expressed in terms of tensor products of $h$, $g$ and $1$, the property of any $R$-form $r$
 $$r(a\otimes bc)=\sum r(\a2\otimes b)r(\a1\otimes c)\quad{\hbox{for every $a,\,b$ and $c\in H_4$}}$$
 implies that it is enough to check that the two forms coincide on $h\otimes h$, $g\otimes g$, $g\otimes h$ and 
 $h\otimes g$.
  Then
 $$\sigma_t\tau*r_s*\sigma^{-1}_t(g\otimes g)=r_{t-s}(g\otimes g)=-1;$$
 $$\sigma_t\tau*r_s*\sigma^{-1}_t(g\otimes h)=\sigma_t\tau*r_s*\sigma^{-1}_t(h\otimes g)=0=
 r_{t-s}(h\otimes g)=r_{t-s}(g\otimes h)$$
  because every summand will involve one of the forms evaluated at a pair composed by the grouplike 
  $g$ and the twisted primitive $h$.  
 $$\sigma_t\tau*r_s*\sigma^{-1}_t(h\otimes h)=
 \sigma_t(h\otimes h)r_s(g\otimes g)\sigma^{-1}_t(g\otimes g)+
 \sigma_t(1\otimes h)r_s(g\otimes h)\sigma^{-1}_t(g\otimes g)+$$
 $$\sigma_t(1\otimes h)r_s(g\otimes 1)\sigma^{-1}_t(g\otimes h)+
 \sigma_t(h\otimes 1)r_s(h\otimes g)\sigma^{-1}_t(g\otimes g)+$$
$$ \sigma_t(1\otimes 1)r_s(h\otimes h)\sigma^{-1}_t(g\otimes g)+
 \sigma_t(1\otimes 1)r_s(h\otimes 1)\sigma^{-1}_t(g\otimes h)+$$
$$ \sigma_t(h\otimes 1)r_s(1\otimes g)\sigma^{-1}_t(h\otimes g)+
 \sigma_t(1\otimes 1)r_s(1\otimes h)\sigma^{-1}_t(h\otimes g)+$$
 $$\sigma_t(1\otimes 1)r_s(1\otimes 1)\sigma^{-1}_t(h\otimes h)=
 -\sigma_t(h\otimes h)+r_s(h\otimes h)+\sigma^{-1}_t(h\otimes h)=$$
 $$-t+s=r_{t-s}(h\otimes h).$$In particular, for $s=t$ one has $\sigma_s\tau*r_s*\sigma_s^{-1}=r_0$.
 
 \begin{remark}{\rm  It is a result by Majid (see \cite{mabook} or \cite{Maji2}) that if two 2-cocycles are cohomologous
 then their corresponding twisted Hopf algebras are isomorphic. If the Hopf algebra involved is dual quasitriangular then the
 corresponding twisted Hopf algebras will
 be isomorphic as dual quasitriangular Hopf algebras. In particular if
 a $2$-cocycle is a coboundary (i.e. it is cohomologous to $\varepsilon\otimes\varepsilon$) one obtains the same dual quasitriangular
 Hopf algebra he started with. The computations above show that the converse  is not true, i.e.
 there are $2$-cocycles which are not coboundaries for which the twist does not change the Hopf algebra structure. For $t\not=0$
 $\sigma_t$ is not a coboundary because otherwise $(H_4,\,r_t)$ would be isomorphic to $(H_4,\,r_0)$  as dual triangular 
 Hopf algebras which is never true by
 the results in \cite{Rad}. }\end{remark}
 We have proved the following
 \begin{theorem}\label{stesso} Let $k$ be a field of characteristic different from $2$.
 Let $H_4$ be Sweedler's Hopf algebra, and $r_t$, $t\in k$ be its universal $R$-forms.
 Then for every $t\in k$, $BC(k,\,H_4,\,r_t)\simeq BC(k,\,H_4,\,r_0)$. The isomorphism maps the class 
 of $A$ in $BC(k,\,H_4,\,r_0)$ to the class of $A_{\sigma_t\tau}$ in $BC(k,\,H_4,\,r_t)$. 
 \hfill$\square$\end{theorem}

 \begin{remark}{\rm If $k$ is algebraically closed, the
 isomorphism between  $BC(k,\,H_4,\,r_s)$ and $BC(k,\,H_4,\,r_t)$ for $st\not=0$ could be deduced
  {\em a priori} by the fact that
 all the dual quasitriangular structures are isomorphic for $st\not=0$, see \cite{Rad}. 
 But $(H_4,\,r_s)\not\simeq(H_4,\,r_0)$ for
  $s\not=0$. Still, we could show that the corresponding Brauer groups are isomorphic. }\end{remark}
  \begin{remark} {\rm Theorem \ref{stesso} combined with the results in \cite{VOZ3} and self-duality of $H_4$
   provide a full descripition
  of $BC(k,H_4,r_t)$ and $BM(k,H_4,R_t)$. This group is the direct sum of $(k,+)$ and the Brauer-Wall group
  of $k$ (which is the Brauer group of the category of ${\bf Z}_2$-graded $k$-modules).
   The computation of the full Brauer group $BQ(k,H_4)$
   and the determination of
  how the different copies of $BC$ fit into $BQ$ in this case is still an open problem.}
  \end{remark}

 \section{The Brauer groups of $H^*$ and $H^{op}$}\label{hop}
 
 In this section we investigate the relation between $BQ(k,\,H)$ and $BQ(k,\,H^*)$ and between
  $BQ(k,\,H)$ and $BQ(k,\,H^{op})$.
  We start with  a variation of a
  Lemma to be found  in \cite{Rad}.
 \begin{lemma} \label{radford} For a finitely generated projective Hopf algebra  
 $H$ over $k$ the usual flip $\tau$ defines an 
 isomorphism between $D(H)$ and $D(H^{op,cop,*})^{op}$. The standard $R$-matrix 
 $\R$ is mapped to $(\tau_{13}\otimes\tau_{24})\R'$
  where $\R'$
 is the standard $R$-matrix of $D(H^{op,cop*})$. In particular, if the antipode $S$ of $H$ is bijective,
  then $\tau(S^{*-1}\otimes S)$ defines an isomorphism between $D(H)$ 
  and $D(H^*)^{op}$ mapping the universal $R$-matrix $\R$ to
  $\tau_{D(H^*),D(H^*)}\R'$.\hfill$\square$\end{lemma}
  This Lemma implies that the categories of left modules of the pairs $(D(H),\,\R)$ and  $(D(H^*)^{op},\,\R')$ are equivalent. 
  Hence the categories of 
  left modules of $D(H)$ and $D(H^*)$ are anti-equivalent braided monoidal categories. 
 Dualizing the above result and using the discussion in the previous sections  and Remark \ref{gio} we have:
  \begin{proposition}\label{maps} Let $H$ be a faithfully projective Hopf algebra with bijective antipode $S$. Then
  $BQ(k,\,H)\simeq BQ(k,\,H^*)$.\end{proposition}
  \Proof $BQ(k,\,H)$ is isomorphic to:
  $$BC(k,\,D(H)^*,\,r)\simeq BC(k,\,D(H^*)^{op*},\,r'\tau)\simeq BC(k,\,D(H^*)^{*cop},\,r'\tau).$$
  Applying the antipode of $D(H^*)^*$ and observing that the $R$-form $r'\tau$ remains unchanged one has:
 $$BQ(k,\,H) \simeq
  BC(k,\,D(H^*)^{*op},\,r'\tau)\simeq BC(k,\,D(H^*)^*,\,r')\simeq BQ(k,\,H^*).$$\hfill$\square$
  
 Observe that the above is in fact an anti-isomorphism. We describe it more explicitely generalizing
 the main
  result in \cite{Til}. 
 We  give it in terms of Yetter-Drinfel'd modules and  
 in terms of modules for the Drinfel'd double.\\
Suppose $(M,\,\cdot,\,\chi)$ is a Yetter-Drinfel'd $H$-module. Then $(M,\,\triangleright)$ is a $D(H)$-module. 
The pull-back of this action along the map $(\tau(S^{*-1}\otimes S))^{-1}=(S^{*}\otimes S^{-1})\tau$ 
defines on $M$ a $D(H^*)^{op}$-module structure that  is given by 
$(a\otimes\varepsilon).m=S^{-1}(a)\cdot m$ for elements of the subalgebra
$H^{op}\otimes \varepsilon$ and 
 by $(1\otimes\xi).m=S^*(\xi)\argh m$ for elements of the sublagebra
$1\otimes H^{*,op}$.
The antipode of $D(H^*)^{op}$ is $S^{-1}\otimes\id$ on $\varepsilon\otimes H^{op}$ and $\id\otimes S^*$ on $1\otimes
 H^{*op}$.
Using the action defined in Remark \ref{gio} the $H^*$-action on $M$ is given by $\cdot$ and $\argh$.
Therefore the map in Proposition \ref{maps} sends a Yetter-Drinfel'd $H$-module 
$(M,\,\cdot,\,\chi)$ to the Yetter-Drinfel'd $H^*$-module $(M,\,\rho,\,\argh)$ where 
 $\rho(m)= \sum m_{(0*)}\otimes m_{(1*)}\in M\otimes H^*$ is  such that for every
  $l\in H\simeq H^{**}$
  one has $(\id\otimes<-,\, l>)\rho(m)=l\cdot m$. This is possible because $M$ is rational.\\
 In other words, the functor from objects in $\ydh$ to objects in $\yddh$ is given by the standard
functors from the category of right $H$-comodules ${\cal M}^H$ to the category of left 
$H^*$-modules $_{H^*}{\cal M}$
and from the category of left $H$-modules  $_H{\cal M}$
  to the category of right $H^*$-modules ${\cal M}^{H^*}$
 associating to a right comodule $(M,\,\chi)$  the left $H^*$-module  $(M,\,\argh)$ 
where $h^*\argh m=(\id\otimes <h^*,\,->)\chi(m)$
and  to the $H$-module $(M,\,\cdot)$ 
  the right $H^*$-comodule $(M,\,\rho)$ where  
 $\rho(m)= \sum m_{(0*)}\otimes m_{(1*)}$ as above. 
One can also check directly that this functor $\cal D$ maps 
Yetter-Drinfel'd modules of $H$ to Yetter-Drinfel'd modules of $H^*$ and that it defines an anti-equivalence of
 braided categories. For instance, the compatibility condition between $\rho$ and $\argh$ for the $H^*$ module 
 $(M,\,\argh)$ and comodule 
 $(M,\,\rho)$ follows from
$$(\id\otimes<-,\, l>)\sum \h1^*\argh \mm0\otimes \h2^*\mm1 $$
$$=\sum \h1^*\argh(\mm0 <\mm1,\,l_{(2)}>) <\h2^*,\,l_{(1)}>$$
$$=\sum \h1^*\argh(l_{(2)}\cdot m) \,<\h2^*,\,l_{(1)}>$$
$$=\sum (l_{(2)}\cdot m)_{(0)}<\h1^*,\,(l_{(2)}\cdot m)_{(1)}>\,< \h2^*,\,l_{(1)}>$$
$$=(\id\otimes <h^*,\,->)\sum (l_{(2)}\cdot m)_{(0)}\otimes (l_{(2)}\cdot m)_{(1)}\,l_{(1)}$$
for any $h^*\in H^*$, $m\in M$ and $l\in H\simeq H^{**}$ together with the
 compatibility condition for the Yetter-Drinfel'd $H$-module $M$.\\
We have seen that 
$\D$ defines a monoidal functor from $\ydh$ to $\yddh^{op}$, the monoidal
 category whose tensor structure is given by $\totimes^{op}$, where $M\tilde\otimes_*^{op}N=N\totimes_* M$
 is the opposite product of $N$ and $M$ viewed as  Yetter-Drinfel'd $H^*$-modules.
The  family of natural isomorphisms 
$$\tau(M,\,N)\colon \D(M)\totimes_*^{op}\D(N)=\D(N)\totimes_*\D(M)\to \D(M\totimes N)$$
compatible with $\D$ and the tensor products 
is given by the usual flips $\tau$. In fact one can also check directly that $\tau$ is a module and comodule isomorphism. For instance
 for $M$  and $N$ objects in  $\ydh$ the left $H^*$-module structure on $M\totimes N$ is\\
 \begin{eqnarray}&&h^*\argh\tau(n\totimes m)=\sum \m0\otimes\n0 <h^*,\,\n1\m1>\\
 &&=\sum \m0 <\h2^*,\,\m1>\otimes \n0 <\h1^*,\,\n1>\\
 &&=\sum (\h2^*\argh_M m)\otimes (\h1^*\argh_N n)
 \end{eqnarray} where $\argh_M$ and $\argh_N$ denote the $H^*$ left module structure associated to the 
 right $H$-comodule structure of $M$ and $N$. The
  right $H^*$-comodule structure $\rho$ of $M\totimes N$  is
 such that for every  $l\in H$ one has
 $$(\id^{\otimes2}\otimes <-,\,l>)\rho(m\totimes n)=l\cdot (m\totimes n)=$$
 $$\sum(\l1\cdot m)\totimes (\l2\cdot n)=\sum (\id\otimes \l1)(\rho_M(m))\otimes 
 (\id\otimes \l2)(\rho_N(n))$$ 
 $$=(\id^{\otimes2}\otimes <-,\,l>)\sum \mm0\otimes\nn0\otimes \mm1\nn1$$ where
 $\rho_M(m)=\sum\mm0\otimes\mm1$ and $\rho_N(n)=\sum\nn0\otimes\nn1$, hence $\tau$ is a comodule isomorphism.\\
One can also see directly that $\D$ and $\tau$ are compatible with  the braidings $\phi$ in $\ydh$ and
$\psi^{op}_{MN}$ in $\yddh^{op}$ where $\psi^{op}_{MN}$ is
$\psi_{NM}\colon N\totimes_*M\to M\totimes_*N$ as Yetter-Drinfel'd $H^*$-modules.
In fact for any pair 
of objects $M$ and $N$ in $\ydh$ and for $n\otimes m\in \D(M)\totimes^{op}\D(N)$,
$$\D(\phi_{MN})\tau_{NM}(n\otimes m)=\phi_{MN}(m\otimes n)=\sum \n0\otimes \n1\cdot m$$
coincides with 
$$\tau_{MN}\circ \psi^{op}_{MN}(n\otimes m)=\tau_{MN}\sum \mm0\otimes <\mm1,\,\n1>\n0=\sum 
\n0\otimes \n1\cdot m.$$
In terms of $\D$ and $\tau$ one sees that if $A$ is a 
Yetter-Drinfel'd $H$-module algebra,
the product in $\D(A)$  is given by $\D(m_A)\tau_{A,A}$, where 
$m_A$ is the product in $A$, i.e. $A$ equipped with the opposite product is a Yetter-Drinfel'd $H^*$-module algebra 
as we had seen before. We have that $\tau$ defines an algebra isomorphism
between $\D(B\# A)=(B\# A)^{op}$ 
and $\D(B)^{op}\#^{op}_*\D(A)\simeq A^{op}\#_* B^{op}$ where $\#_*$ the tensor product 
of Yetter-Drinfel'd $H^*$-module algebras. This can be checked directly 
 for every $a,\,c\in A$ and $b,\,d\in B$: 
 $$m_{(B\# A)^{op}}(\tau\otimes\tau)((a\#_*b)\otimes(c\#_*d))= d\phi_{AB}(c\otimes b)a=$$
$$  \tau\bigl(\sum a\circ \tau\phi_{AB}(c\otimes b)\circ d\bigr)= 
\tau\bigl(a\circ(\psi_{BA}\tau(c\otimes b))\circ d\bigr)$$
 $$=\tau\bigl(a\circ(\psi_{BA}(b\otimes c))\circ d\bigr)=\tau m_{A^{op}\#_* B^{op}}((a\#_*b)\otimes(c\#_* d))$$
  where $\circ$ denotes the opposite product in $A$ and $B$ and the interchange between $\tau\phi_{AB}$ and 
 $\phi_{BA}\tau$ follows by the compatibility of $\D$ and $\tau$ with the braidings.
\begin{proposition}\label{isomorphic} Let $H$ be a finitely generated  projective Hopf algebra over $k$ 
with bijective antipode. Then the isomorphism between $BQ(k,\,H)$ and $BQ(k,\,H^*)$ of Proposition \ref{maps}
maps $[A]\mapsto[\overline{\D(A)}]_*$. 
Here $[B]_*$ denotes the equivalence class of $B$ as $H^*$-Azumaya algebra.
\hfill$\square$
\end{proposition}
For sake of completeness we show also  the (functorial) reciprocity for the opposite algebra,
 $F$ and $G$ with respect to $\D$. 
$\overline{\D(A)}$ has product $m_{\D(A)}\psi_{AA}^{op}=m_A\tau\psi_{AA}=m_A\phi_{AA}\tau=m_{\D(\bar A)}$.
  Hence $\D(\bar A)\simeq\overline{\D(A)}$.
As far as the canonical maps in $\yddh$ $F_*$ and $G_*$ are concerned, $F_*$ is the composite:
$$\D(A)\#^{op}_*\overline{\D(A)}=\D(A)\#_*\D(\bar A)\buildrel\tau \over\longrightarrow\D(A\#\bar A)\buildrel
\D (F)\over\longrightarrow\D(End(A))=End(\D(A))^{op}$$ and similarly for $G_*$
$$\overline{\D(A)}\#^{op}_*\D(A)=\D(\bar A)\#_*\D(A)\buildrel\tau \over\longrightarrow\D(\bar A\# A)\buildrel
\D (G)\over\longrightarrow\D(End(A)^{op})=End(\D(A)).$$
It is a small exercise in sigma notation to check that $F_*(a\#_*b)(c)=G(\bar b\# a)(c)$,
 $G_*(\bar b\#_* a)(c)=
F(a\# \bar b)(c)$
and that $\D(End(A))$ is indeed $End(A)^{op}$  with the natural Yetter-Drinfel'd $H^*$-module structure.  

\begin{corollary} For a finitely generated projective Hopf algebra $H$ over a commutative ring $k$, 
with bijective
 antipode the isomorphism in Proposition \ref{isomorphic} induces an isomorphism between 
$BQS(k,\,H)$ and $BQS(k,\,H^*)$.
\end{corollary}
\Proof It is immediate to check that $k$ is a direct summand of $A$ as a Yetter-Drinfel'd $H$-module if and only if 
$k$
is a direct summand of $A$ as a Yetter-Drinfel'd $H^*$-module.\hfill$\square$

\begin{corollary} If $H$ is commutative {\em or} cocommutative then all the ``associated'' Hopf algebras
 (opposite, co-opposite, 
opposite of the dual, co-opposite of the dual, etcetera) have the same Brauer group.
\hfill$\square$
\end{corollary}
It is an interesting question whether an isomorphism between  $BQ(k,\,H)$ and $BQ(k,\,H^{op})$ holds in general. 
Here follow some observations concerning this question.
\begin{lemma}\label{doublecop} The linear map
  $$S^*_H\otimes\id\colon D(H^{cop})\to D(H)^{cop}$$ is a Hopf algebra isomorphism. If for $D(H)^{cop}$ 
  we have the 
  $R$-matrix  $\tau \R_{D(H)}$ then the corresponding $R$-matrix in $D(H^{cop})$ is 
  $\tau \R^{-1}_{D(H^{cop})}$.\end{lemma}
\Proof Direct computation.\hfill$\square$
  
 \begin{proposition}\label{cop} For any finitely generated projective Hopf algebra $H$, there is an isomorphism 
  between $BQ(k,\,H)$ and $BM(k,\,D(H^{cop}),\,\tau \R^{-1}_{D(H^{cop})})$.\end{proposition}
  \Proof One has
   $$BQ(k,\,H)\simeq BM(k,\,D(H),\,R)\simeq BM(k,\,D(H)^{cop},\,\tau \R)$$ where the second
    isomorphism is clear because
  the algebra structure is unchanged. By Lemma \ref{doublecop}
  $$BQ(k,\,H)\simeq BM(k,D(H^{cop}),\tau\R^{-1}_{D(H^{cop})})\simeq 
  BM(k,D(H^{op}),\tau \R^{-1}_{D(H^{op})})$$ where the second isomorphism is induced by  
  the isomorphism $S_H^{-1}\otimes S^*_H$ 
  from $D(H^{op})$ to
  $D(H^{cop})$ that leaves the $R$-matrix invariant.\hfill$\square$
  
  Under the above series of maps (coming from equivalences of categories) the Yetter-Drinfel'd
   $H$-module $(M,\,\cdot,\,\chi)$ 
  is  mapped to the Yetter-Drinfel'd $H^{op}$-module $(M,\,\rightharpoondown,\,\chi)$ with 
   $h\rightharpoondown m:=S(h)\cdot m$  where $S$ is the antipode of $H$.
   In fact 
 $$\chi(h\rightharpoondown m)=\chi(Sh\cdot m)=\sum (Sh)_{(2)}\cdot\m0\otimes (Sh)_{(3)}\m1 S^{-1}(Sh)_{(1)}=$$
 $$\sum S(h)_{(2)}\cdot\m0\otimes S(h_{(1)})\m1 h_{(3)}=\sum S(h)_{(2)}\cdot\m0\otimes  h_{(3)}\circ\m1
\circ S^{-1}_{H^{op}}(h_{(1)})$$ where $\circ$ denotes the opposite product. 
This construction defines a functor which, together with the natural transformation $\tau$
  gives an anti-equivalence of monoidal categories. The natural isomorphism $\tau$ would respect the braiding
   if for every Yetter-Drinfel'd module
   $M$ and $N$, for every $m\in M$ and $n\in N$ one had
  $$\sum\m0\otimes\m1\cdot n=\sum (S\n1)\cdot m\otimes \n0.$$ 
 In terms of  modules over Drinfel'd quantum double of $H$,
  $$\sum\m0\otimes\m1\cdot n=\tau(\R\triangleright(n\otimes m))$$
  if $\R=\sum \R_1\otimes \R_2$ is the canonical $\R$-matrix for $D(H)$. If
  $P=\sum P_1\otimes P_2$ is
  $\tau_{D(H),D(H)}\R^{-1}=\tau_{D(H),D(H)}(\id\otimes S^{-1}_{D(H)})\R$, then 
 $$\sum (S\n1)\cdot m\otimes \n0=\tau(P\triangleright(n\otimes m)).$$
  Hence $\tau$ would respect the braidings if $D(H)$ were triangular. Then  
  $BQ(k,\,H)$ would be isomorphic to $BQ(k,\,H^{op})$. 
 However, for $k$ a field and $H$ nontrivial, $(D(H),\R)$ is never triangular.
  An easy way to see this is given by the theory of the exponent of a 
  finite-dimensional Hopf algebra introduced by 
  P. Etingof and S. Gelaki
  in \cite{geleti}. The authors show  the so-called exponent of a Hopf algebra $H$ is equal to
    the order of $(\tau \R)\R$ in $D(H)^{\otimes 2}$
   and also equal to the order of the element
    $u=m(S_{D(H)}\otimes \id^{\otimes 2})\tau(\R)=\sum S^{*-1}(h_j^*)\otimes h_j$ in $D(H)$.
  Here  
 $\{h^*_j\}_{j\in J}$ and $\{h_j\}_{j\in J}$ are dual bases and $S^*$ is the antipode of $H^*$.
  It is clear then that $(\tau \R)\R\not=\varepsilon\otimes 1\otimes \varepsilon\otimes 1$ because 
  $u\not=\varepsilon\otimes 1$.

 Observe that triangularity of $\R$ is not a necessary condition for $BQ(k,\,H)\simeq BQ(k,\,H^{op})$ to hold.
  This is clear by looking at the case of  $H$ commutative. In particular Proposition \ref{cop} implies that if $H$ is commutative then
  there is an isomorphism between $BC(k, D(H), \R)$ and $BC(k, D(H), \tau \R^{-1})$.\vskip20pt
 

\centerline{\bf ACKNOWLEDGEMENTS}\vskip7pt 

This work has been carried out at the University of Antwerp where
 the author had a post-doc
position 
financed by the EC network ``Algebraic Lie
Representations''contract ERB-FMRX-CT97-0100.
The author wishes to thank University of Antwerpen for the warm hospitality and  Professor Fred
 Van Oystaeyen, Dr. Juan Cuadra
and Dr. Yinhuo Zhang for valuable comments and useful discussions.

\end{document}